\documentclass{article}
\usepackage{arxiv}

\usepackage[utf8]{inputenc} 
\usepackage[T1]{fontenc}    
\usepackage{hyperref}       
\usepackage{url}            
\usepackage{booktabs}       
\usepackage{amsfonts}       
\usepackage{nicefrac}       
\usepackage{microtype}      
\usepackage{lipsum}		
\usepackage{graphicx}
\usepackage{natbib}
\usepackage{doi}

\usepackage{float}
\restylefloat{table}

\usepackage{parskip}
\usepackage{setspace}
\newcommand{\Tr}{\operatorname{Tr}}
\newcommand{\Det}{\operatorname{Det}}
\usepackage{amsthm,amsmath,color,mathtools,enumerate}
\usepackage{amsfonts,amssymb,graphicx,float}
\usepackage[math]{iwona}
\usepackage{mathtools}
\DeclarePairedDelimiter{\norm}{\lVert}{\rVert}
\usepackage{eufrak,amscd,bezier,latexsym,mathrsfs,enumerate,multirow}
\usepackage{enumitem}
\newtheorem{theorem}{Theorem}
\newtheorem{corollary}[theorem]{Corollary}
\title{A comprehensive numerical investigation of a coupled mathematical model of neuronal excitability}

\author{ \href{https://orcid.org/0000-0002-4253-5877}{\includegraphics[scale=0.06]{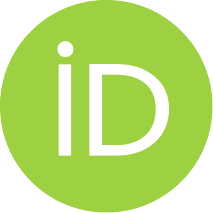}\hspace{1mm}Burcu~G\"urb\"uz}\thanks{Corresponding author}\\
	Institute of Mathematics\\
	Johannes Gutenberg-University Mainz\\
	55128 Mainz, Germany \\~\\
	Institute for Quantitative and Computational Biosciences (IQCB)\\
	Johannes Gutenberg-University Mainz\\
	55128 Mainz, Germany \\
	\texttt{burcu.gurbuz@uni-mainz.de} \\
	\And
	\href{https://orcid.org/0000-0003-1421-3966}{\includegraphics[scale=0.06]{orcid.pdf}\hspace{1mm}Ayt\"{u}l~G\"{o}k\c{c}e} \\
	Department of Mathematics\\
	Ordu University\\
	52200 Ordu, Turkey \\
	\texttt{aytulgokce@odu.edu.tr} \\
	\And
	\href{https://orcid.org/0000-0002-7743-3512}{\includegraphics[scale=0.06]{orcid.pdf}\hspace{1mm}Mahmut~Modanl{\i}} \\
	Department of Mathematics\\
	Harran University\\
	63300 \c{S}anl{\i}urfa, Turkey \\
	\texttt{mmodanli@harran.edu.tr} \\
}

\renewcommand{\shorttitle}{\textit{arXiv} Template}

\begin{document}
\maketitle

\begin{abstract}
Being an example for a relaxation oscillator,  the FitzHugh-Nagumo model has been widely employed for describing the generation of action potentials. In this paper, we begin with a biological interpretation of what the subsequent mathematical and numerical analyses of the model entail. The interaction between action potential variable and recovery variable is then revisited through linear stability analysis around the equilibrium and local stability conditions are determined. Analytical results are compared with numerical simulations. The study aims to show an alternative approach regarding Taylor polynomials and constructed difference scheme which play a key role in the numerical approach for the problem. The robustness of the schemes is investigated in terms of convergency and stability of the techniques. This systematic approach by the combination of numerical techniques provides beneficial results which are uniquely designed for the FitzHugh-Nagumo model. We describe the matrix representations with the collocation points. Then the method is applied in order to acquire a system of nonlinear algebraic equations. On the other hand, we apply finite difference scheme and its stability is also performed. Moreover, the numerical simulations are shown. Consequently, a comprehensive investigation of the related model is examined.
\end{abstract}



\section{Introduction}
\label{Sec:Intro}

Tracing back its roots to seminal works by Richard FitzHugh in $1961$ \cite{fitzhugh1961impulses} and Jin-Ichi Nagumo in $1962$ \cite{nagumo1962active}, the FitzHugh-Nagumo model (FHN) is often considered to model excitable systems and represents a two component simplified version of well known Hodgkin-Huxley model (HH) \cite{rocsoreanu2012fitzhugh,tuckwell1998analytical,jones1984stability,sato1992response}. This transition to two variable model developed by FitzHugh  was made possible through The van der Pol  equations to understand nonlinear relaxation oscillators \cite{fitzhugh1961impulses,izhikevich2006fitzhugh}. Here, a further analysis of this model with the addition of several terms inspired FitzHugh to expand this knowledge with  a simplified adaptation of the Bonhoeffer van der Pol model (BVP), that has been renowned as FitzHugh-Nagumo model \cite{fitzhugh1961impulses,izhikevich2006fitzhugh}. Although the development of biophysically more realistic HH model was a breakthrough to understand how action potentials arise and propagate, two variable FHN model  provided simple and conceptual approach  that would not only capture the fundamental aspects surrounding neuronal excitability, but also allow for phase plane methods to gain a better understanding of spike generation \cite{tuckwell1998analytical,sherwood2014fitzhugh,faghih2010fitzhugh,kostova2004fitzhugh}.

Although the FHN model is originally use to understand the action potential of single neurons, it has been broadly exploited for many reaction diffusion systems whose components comprise a fast activator variable  and a slow inhibitor variable. Thus, the model can be amenable to deal with various dynamics in excitable media and still described as a corner stone for models that exhibit similar characteristics \cite{onyejekwe2021dynamical,izhikevich2003bursts}. In addition to capturing neural spiking behaviours, this model has been adapted to describe a wide spectrum of
excitable media from the Belousov-Zhabotinsky chemical reaction \cite{pertsov1984rotating} to slime mold amoeba \cite{vasiev1994simulation} to cardiac tissue \cite{aliev1996simple,glass1995resetting,panfilov1995re}.

The generic FHN model is driven by two variables of state, representing both membrane potential and refractoriness (or recovery). Equipped with the FHN model, we investigate its stability by considering a small perturbation around the equilibrium point and its corresponding Taylor expansion. By expressing the linearised system in vector form with the Jacobian matrix evaluated at this point, we seek exponential solutions; this leads to the characteristic polynomial. This approach is followed by many researchers to gain better insights about solution states of real life problems, see for example \cite{murray2001mathematical,liao2007stability,seydel2009practical} for a comprehensive overview. Due to the complexity of the characteristic equation, we firstly consider that no input current is applied. Once this is achieved, we derive conditions needed for various dynamical characteristics and explore the model in more depth for which we conduct bifurcation analysis when non-zero current is applied.
We expand upon this extensive body of knowledge by refining the FHN model further.

The numerical investigation of the FHN model is also important to understand the dynamical properties of the model. The unknown functions play a crucial role as well as the parameters for a comprehensive approach of the model. Thus, the applications of the model with the help of numerical approaches show the results of the unknowns in a particular time interval. On the other hand, the model has been established in different types of differential equations and systems. Hence the details of the model has been considered with the different aspects. The simplified version of the model gives us understanding in the application of neural networks in biological neural concept, cardiovascular systems, deep learning and so on. The numerical investigation of such systems have been investigated by several authors including \cite{olmos09,chap,feng,q1,q2}. The fractional FHN model has been solved numerically by an implicit numerical method with the shifted Grunwald-Letnikov approximation \cite{liu}, Jacobi-Gauss-Lobatto collocation method has been used for  the solution of the generalized FHN equation together with time-dependent coefficients \cite{bra}, for solving the stochastic FHN systems some numerical and analytical methods have been improved by the authors \cite{li,tuck}, also traveling wave solutions of FitzHugh model with cross-diffusion \cite{bm6} and multifront regime of a piecewise-linear FHN model with cross diffusion has been presented \cite{bm7}. Besides, bifurcation analysis of a FHN model \cite{mm1}, hidden extreme multistability and synchronicity of memristor-coupled non-autonomous memristive FHN models \cite{mm2} have been investigated. The main idea of the applications of these numerical approaches to obtain accurate results with the help of straightforward algorithmic steps based on the present model.
On the other hand, numerical techniques have been implemented in various mathematical model investigations such as a trigonometric quintic B-spline method has been presented for the solution of a class of turning point singularly perturbed boundary value problems (SP-BVPs) whose solution exhibits either twin boundary layers near both endpoints of the interval under consideration or an interior layer near the turning point \cite{mm3}. Recently, many studies have been carried out on numerical solutions of the mathematical problems by using different methods \cite{mm4, reviewers1,reviewers2,reviewers3}. In our research, we explore an alternative numerical approach that combines a Taylor series method with a specially constructed difference scheme. This method provides highly accurate, well established, convergent and stable results under reasonable computational conditions for solving the FHN model.

One popular relaxation oscillator for explaining how action potentials are generated in biological systems is the FHN model. Still, a thorough investigation of both analytical and numerical methods is still needed. By re-examining the relationship between the action potential and the recovery variables through linear stability analysis, this study aims to fill this gap by providing a biological explanation to support the mathematical and numerical analyses.
In particular, the work presents an alternative numerical method specifically tailored for the FHN model, which makes use of the strength of Taylor polynomials and developed difference schemes. By combining these analytical and numerical methods, the work aims to give a thorough grasp of the behaviour of the model, providing information on its stability, convergence, and robustness properties. In the numerical technique, the matrix representations, collocation points, and resulting system of nonlinear algebraic equations are presented in detail. The stability of the finite difference scheme is also investigated, which further strengthens the resilience of the numerical approach. Finally, this study aims to add to the existing body of knowledge on relaxation oscillators and action potential generation by offering a distinct and methodical analysis of the FHN model using a complementary blend of analytical and numerical methods \cite{yonet23,gurbuzg22, modanli23}.

The organisation of the paper is the following. In Section \ref{model}, the FHN model is revisited and preliminary mathematical results are presented. Here, the stability around the possible equilibria is determined via linearisation method. Eigenvalues are computed and dynamical behaviour is shown in the absence and presence of external input. Then, in Section \ref{Sec:Taylor}, Taylor expansion to the FHN model is presented with the error estimation for the numerical solution and the convergence of the numerical technique based on Taylor truncated series are discussed. Next, in Section \ref{Sec:fdifference}, the Difference scheme approach is constructed as an alternative route for stability analysis of the system. Finally in Section \ref{Sec:Conc}, summary of the results are given and potential future directions are discussed.

\section{Revisiting the FHN model: preliminary mathematical  results}\label{model}
As mentioned in Section \ref{Sec:Intro}, the two-variable FHN model
provides a simplified but analytically-robust platform for extrapolating the central features
of the HH model  \cite{fitzhugh1961impulses,izhikevich2006fitzhugh,faghih2010fitzhugh,kostova2004fitzhugh}. The commonly used form of FHN model
is given by
\begin{equation}\label{eq-model1}
\setstretch{1.4}\begin{cases}
\mu \dfrac{dv}{dt}=f(v,a)-\omega+I \equiv \mu  F(v,\omega),\\
\dfrac{d\omega}{dt}=v-\gamma \omega \equiv G(v,\omega),
\end{cases}
\end{equation}
where $f(v,a)=v(a-v)(v-1)$  with $0<a<1$, $\mu,\gamma>0$.

\begin{table}[h]
\caption{ Parameters and variables used in the model \eqref{eq-model1}.
\label{Tab:1}}
\begin{center}
\begin{tabular}{ c l}
\hline
parameter/variables  & biological meaning \\
\hline
$v$ &  membrane potential  \\
$\omega$ &  recovery (adaptation) variable \\
$I$ &  input current\\
$\mu$ & time course for membrane potential   \\
$\gamma$ & depletion strength of recovery variable \\
$a$ &  threshold between electrical silence and
electrical firing \\
\hline
\end{tabular}
\end{center}
\end{table}
Now, we briefly revisit the equilibrium, stability and bifurcation studies of the FHN model, which has been extensively studied analytically and numerically in literature.
\subsection{Equilibrium and stability of the model with fixed parameters}
The equilibrium point  $E_\ast= E(v_\ast,\omega_\ast) $ of the system is found using
\begin{equation}
   f(v_\ast,a) - \omega_\ast+I =0\;\;\textrm{and}\;\; v_\ast- \gamma \omega_\ast=0,
\end{equation}
from which
\begin{equation}
f(v_\ast,a)- \dfrac{v_\ast}{\gamma} +I =0,
\label{Equ:steady}
\end{equation}
is obtained.
One can draw conclusions about the stability at this point by regarding a small perturbation
of the form
\begin{equation}
(v,\omega) = (v_\ast,\omega_\ast) + \epsilon (\psi_1(t), \psi_2(t)) + O(\epsilon^2), \;\; 0<\epsilon<< 1.
\end{equation}
Using a first order Taylor expansion
and substituting in the model \eqref{eq-model1} the system is presented as
\begin{equation}
\begin{aligned}
\frac{d \psi_1}{d t} & = \frac{f'(v_\ast,a) \psi_1(t) -\psi_2(t)}{\mu},\\
\frac{d \psi_2}{d t} &= \psi_1(t)- \gamma \psi_2(t),
\end{aligned}
\label{Eq:Linearisation}
\end{equation}
where $f'$ represents the derivative of function with respect to $v$. The matrix form of the system in \eqref{Eq:Linearisation} is written as
\begin{equation}
 \dfrac{d}{d t}\begin{pmatrix}
   \psi_1  \\
   \psi_2 \\
     \end{pmatrix}=
     \left. \mathcal{M} \right\rvert_{E_\ast}
     \begin{pmatrix}
   \psi_1  \\
   \psi_2 \\
     \end{pmatrix}, \;\;\;\left. \mathcal{M}  \right\rvert_{E_\ast}=  \begin{pmatrix}
   \dfrac{f'(v_\ast,a)}{\mu} & \dfrac{1}{\mu} \\
   1 & -\gamma \\
     \end{pmatrix},
\end{equation}
where $\mathcal{M}$ represents the Jacobian matrix evaluated at equilibrium $(v_{\ast},\omega_{\ast})$. To find non-trivial solutions,  the characteristic polynomial $\Det (\left. \mathcal{M}  \right\rvert_{E_\ast}- \lambda I) =0$ is solved for which the eigenvalues satisfy
\begin{equation}
 2 \lambda_{1,2} = - \Tr \left.\mathcal{M}  \right\rvert_{E_\ast} \pm \sqrt{\Tr \left. \mathcal{M} ^2 \right\rvert_{E_\ast}  - 4 \Det \left. \mathcal{M}  \right\rvert_{E_\ast}} .
 \label{Equ: Eigenvalue}
\end{equation}
Note that, in the absence of synaptic current, i.e. $I=0$, we have
\begin{equation}
f(v_\ast,a) - \dfrac{v_\ast}{\gamma} = - v_\ast \left\{v_\ast^2 - (a+1) v_\ast + \left(a+ \dfrac{1}{\gamma}\right)\right\}=0,
\end{equation}
leading to $v_\ast^{(1)} =0$ and
\[v_\ast^{(2,3)} = \dfrac{(a+1)}{2} \pm \dfrac{\sqrt{(1-a)^2 - 4/\gamma }}{2}.\]
This implies that real $ v_\ast^{(2)}$ and $v_\ast^{(3)}$ exist only if  $\gamma (1-a)^2 \geq 4$. Thus $a$ and $\gamma$ parameters are indicative of the existence of non-zero membrane potential. The case  $\gamma (1-a)^2 <4$ corresponds to extinction state. The Jacobian matrix for the extinction state
\begin{equation}
\left. \mathcal{M}  \right\rvert_{(0,0)}  = \begin{pmatrix}
\dfrac{-a}{\mu} & \dfrac{-1}{\mu}\\
1 & -\gamma
\end{pmatrix},
\end{equation}
implies that  $\Tr \left. \mathcal{M}  \right\rvert_{E_0} <0$ and  $\Det \left. \mathcal{M}  \right\rvert_{E_0} >0$ for which $E_0=(0, 0)$. Thus the extinct  state  is either a stable node or spiral depending on the value inside the square root in equation \eqref{Equ: Eigenvalue}.
In the case where $I\neq0$, it is possible to have  $\Tr \mathcal{M}>0$ since the equilibrium of the system introduced by  equation \eqref{Equ:steady} depends on the input current $I$.  We  may now look for the Hopf bifurcations that may arise when $\Tr \mathcal{M} =0$. Since the equilibrium is stable when $ I = 0$, we expect an instability as $I$ increases, followed by a
restabilization for a larger value of $I$. It is worth to note that the existence of the equilibria $v_\ast^{(2,3)}$ depend on $a$ and $\gamma$ parameters and the gradient of the $w$-nullcline decreases as $\gamma$ is increased. For  $I = 0$, the nullclines intersect at
only singular point (a resting state), which is a stable node. Thus, if $(1-a)^2<4/\gamma$ then we only have one equilibrium point at $(v_\ast,\omega_\ast)=(0,0)$ and $f'(0)=-a$.

To validate the above analytical study, the numerical simulations of the FHN system given in \eqref{eq-model1} can be shown. Throughout the rest of the paper we keep the system parameters constant at $a=0.22$, $\gamma=1.18$ and $\mu=0.008$, while varying the input current $I$. The initial conditions (ICs) for the system are initialised to $(v(0), \omega(0)) = (0, -0.2)$. Besides, the initial conditions (ICs) for the system are set to $(v(0),\omega(0))=(0,-0.2)$.

\begin{figure}[ht!]
\centering
\includegraphics[width=1\textwidth]{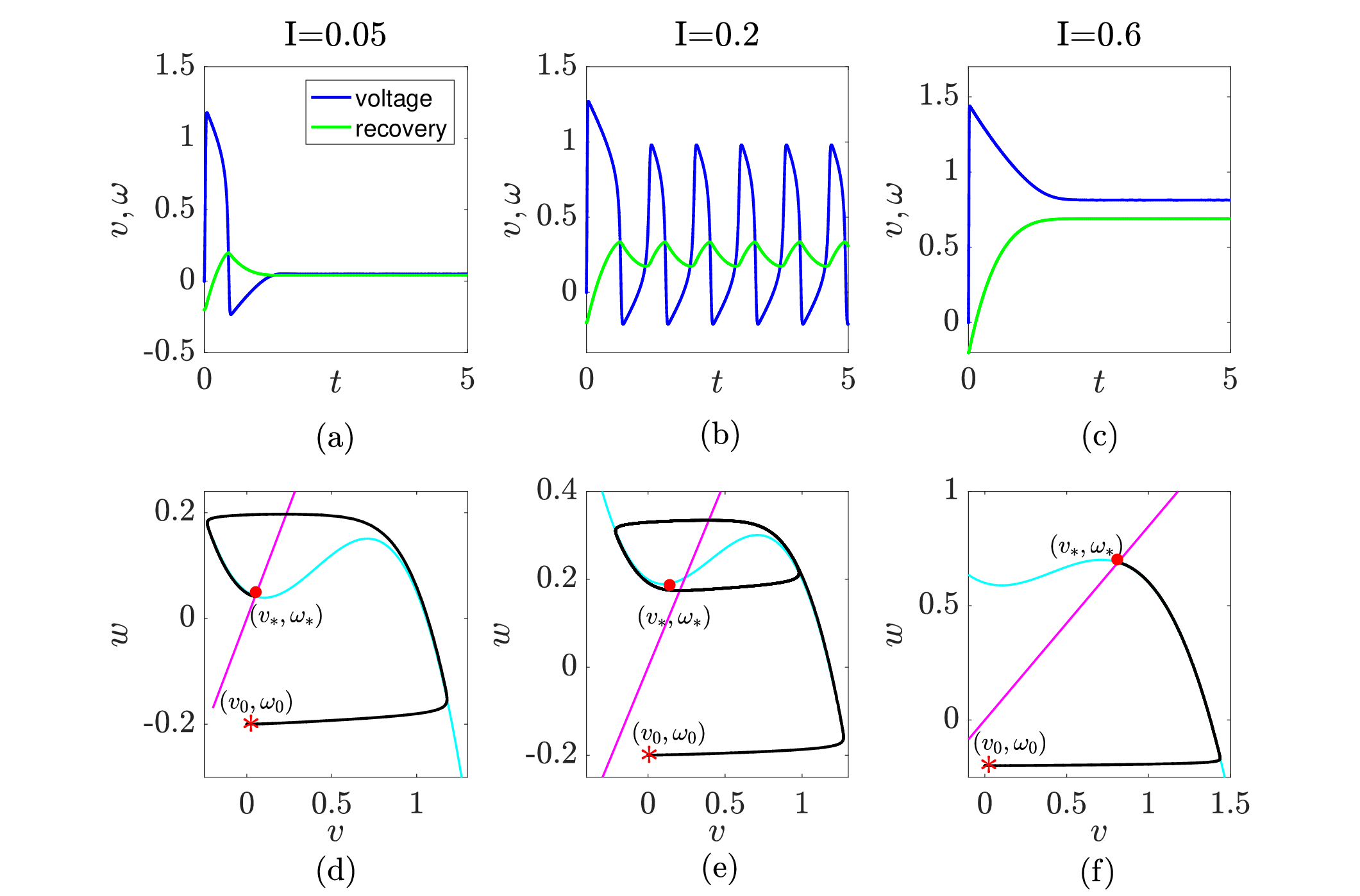}
\caption{ The trajectory of the voltage variable (blue) and recovery (green) variable are shown for different levels of input, where $I=0.05$ (a,d), $I=0.2$ (b,e) and $I=0.6$ (c,f). Here (d-f) represent the corresponding phase diagrams. The red star represents the IC, that is chosen as $(v_0,\omega_0)=(0,-0.2)$. The equilibrium points for each case are shown with red dot.  The dashed magenta and cyan curves represent nullclines of the system other parameters are given in the text.
}
\label{Fig:TimeEvolution}
\end{figure}

Figure \ref{Fig:TimeEvolution} represents recovery variable $\omega$ (a-c) and the time evolution of membrane potential $v$ and the corresponding phase-plane diagrams (d-f) with different levels of input currents: $I=0.05$ (a,d), $I=0.2$ (b,e) and $I=0.6$ (c,f). Here Fig. \ref{Fig:TimeEvolution}(a,d) corresponds to excitable system where the fixed point $(v_\ast,\omega_\ast)=(0.0495,0.042)$ is a stable spiral. Figure \ref{Fig:TimeEvolution}(b,e) show  the oscillatory behaviour of the system when input $I$ crosses the first Hopf bifurcation. Here the equilibrium point is found to be $(v_\ast,\omega_\ast)=(0.2404,0.2037)$ with $I=0.2$.   In addition, one can observe how with large current input the oscillations drop to the steady state $(v_\ast,\omega_\ast)=(0.8141,0.6899)$ after passing the other bifurcation point  in   Fig. \ref{Fig:TimeEvolution}(c,f). The dashed magenta and the nullclines of the system are represented by cyan curves.

\begin{figure}[ht!]
\centering
\includegraphics[width=0.8\textwidth]{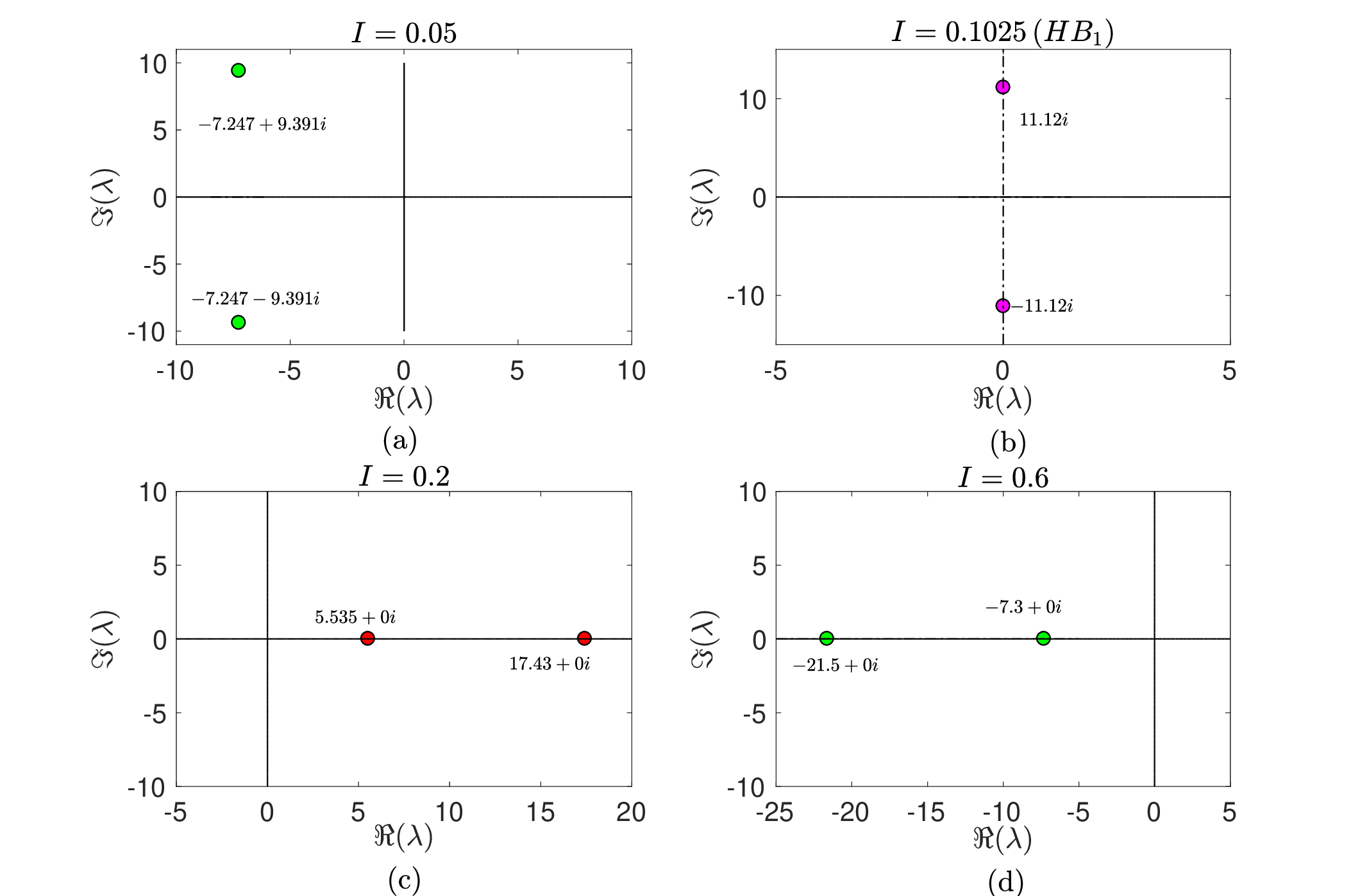}
\caption{Spectra of steady state $(v_\ast,\omega_\ast)$ for FHN system given in \eqref{eq-model1} for different values of input $I$. The system is stable spiral when $I=0.05$ (a), The first Hopf bifurcation $(HB_1)$ occurs at $I=0.1025$ (b), where a pair of eigenvalues with only imaginary part is obtained. The second Hopf bifurcation ($HB_2)$ occurs at $I=0.4963$ (not shown). The system represents an unstable node for  $I=0.2$ (c) and stable node for $I=0.6$(d).}
\label{Fig:Eigens}
\end{figure}

Figure \ref{Fig:Eigens} presents the roots of the characteristic equation, given by equation \eqref{Equ: Eigenvalue} for varying values of $I$ parameter. When a small input current is incorporated in the model, two eigenvalues have negative real parts (green dots) therefore the system is a stable spiral, see Fig. \ref{Fig:Eigens}(a). Then the first Hopf bifurcation ($HB_1$) occurs at a critical threshold value  $I=0.1025$, as seen in  Fig. \ref{Fig:Eigens}(a), where only a pair of imaginary ($\lambda= \pm 11.12i$) eigenvalues (magenta dots) is observed . After crossing $HB_1$, e.g. for  $I=0.2$ in Figure \ref{Fig:Eigens}(c), the positive coexisting state becomes an unstable node  until it passes the second Hopf bifurcation ($HB_2$), that occurs at $I=0.4963$ (not shown). Increasing the input current further, e.g. $I=0.6$, the stability switches back from unstable to stable, where the coexisting state becomes a stable node.

\subsection{Bifurcation analysis}
The analytical and numerical studies presented in this section of the model have been validated by numerous works published in the literature over many years, see for example \cite{sherwood2014encyclopedia, cebrian2024six}.

For the case where  $I\neq0$, $\Tr \mathcal{M}$ could be positive as the equilibrium point given  by equation \eqref{Equ:steady} depends on the input current. Since the equilibrium point is stable when $I=0$, we expect an instability when input current $I$ increases, followed by a restabilisation for larger values of $I$, as stated in Figs. \ref{Fig:TimeEvolution} and \ref{Fig:Eigens}.

\begin{figure}[ht!]
\centering
\includegraphics[width=1\textwidth]{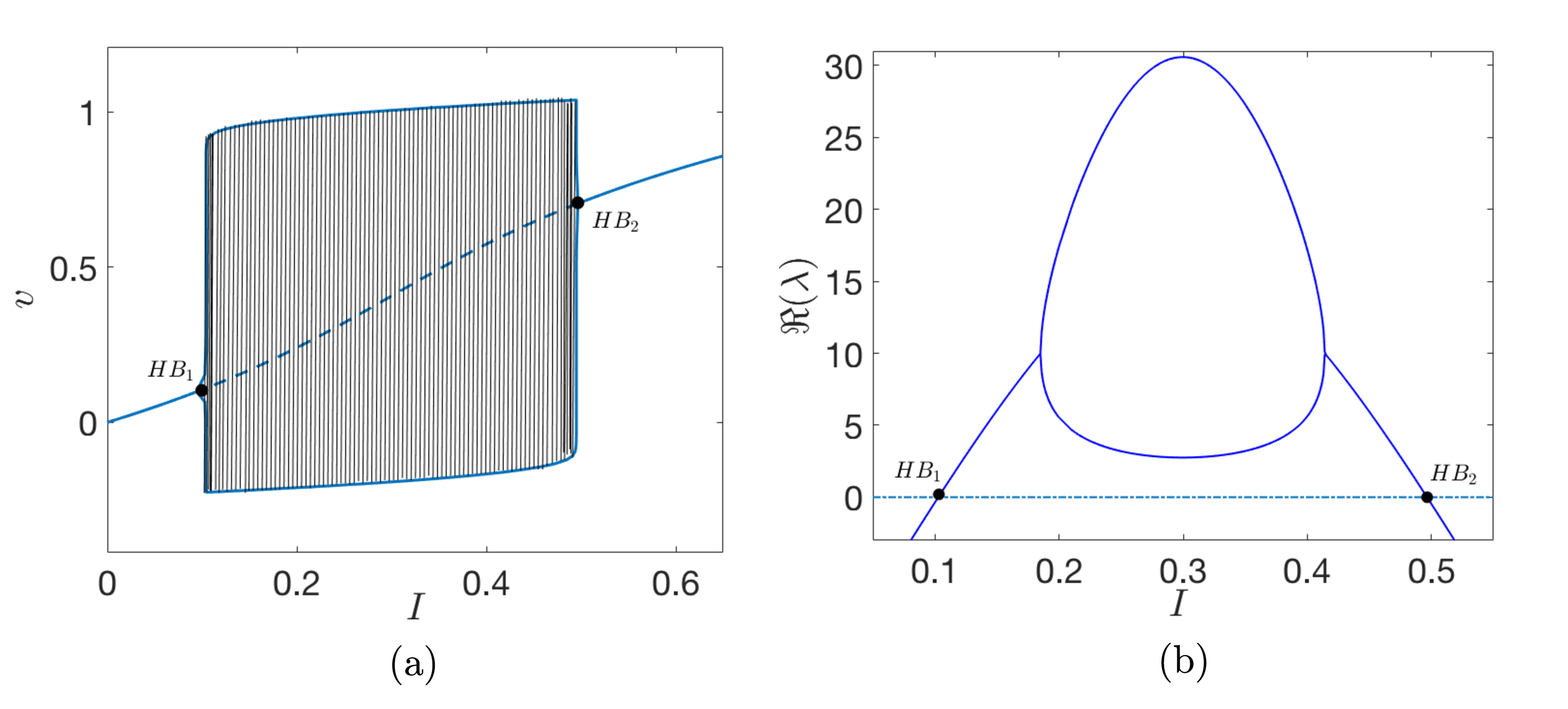}
\caption{Bifurcation diagram of $v$ with respect to input $I$ (a), where solid line show stable dynamics and dashed line shows unstable dynamics with two eigenvalues with positive real parts and periodic orbits arising from $HB_1$ and $HB_2$ are observed. The real parts of the roots of the characteristic equation are plotted with respect to input $I$ (b). }
\label{Fig:Bifurcation}
\end{figure}

Bifurcation diagram of the system \eqref{eq-model1} is shown in Fig. \ref{Fig:Bifurcation}(a). Here the system may be considered to be excitable when the coexisting state lies on the left  end of the interval for  input current $I$ (stable). On the contrary, it becomes  oscillatory when the equilibrium  point lies on the middle branch with dashed line (unstable). The equilibrium looses its stability through the first Hopf bifurcation ($HB_1$)  at $ \lambda_{1,2} = \pm  11.12 i $, leading to a stable limit cycle surrounding unstable fixed point. Increasing the values of input current from $I=0.1025$ to $I=0.4963$, a small increase occurs in size of the limit cycle for higher values of input current and the system evolves to a limit cycle with a corresponding change in its size.
The second Hopf bifurcation ($HB_2$) is encountered  when input reaches to $I=0.4963$, by which the oscillations drop to a coexisting equilibrium. Note that both bifurcations are supercritical giving rise to stable limit cycle. This can be also validated depending on the sign of the first Lyapunov coefficient \cite{engelborghs2002numerical}. These results are also confirmed  in Fig. \ref{Fig:Bifurcation}(b), where real part of the eigenvalues are plotted with respect to input parameter and stability switch is observed at  two critical thresholds for Hopf points.

The bifurcation diagrams in Fig. \ref{Fig:Bifurcation} are plotted using DDE-BIFTOOL package, comprising a collection of MATLAB functions and allowing the stability analysis of fixed points and the numerical continuation of the model around the fixed points \cite{engelborghs2002numerical}. This package is commonly used to study continuation of the delay differential equations, yet it can also be used for systems with no delay terms. The FHN model is a mathematical model used to describe the dynamics of excitable cells, particularly neurons. The model describes the behaviour of the membrane potential ($v$) and a recovery variable ($\omega$). Biologically, the FHN model captures the essential dynamics of neuronal firing. When $v$ exceeds a certain threshold, it triggers an action potential (excitation) leading to a rapid depolarisation phase. Subsequently, $\omega$ influences the recovery phase, causing the membrane potential to return to its resting state (inhibition). Therefore, the dynamic behaviour of the model has been obtained analytically and simulated in the figures above. This interplay between $v$ and $\omega$ is similar to the behaviour of neurons during the firing and recovery phases of action potentials. Overall, the FHN model provides a simplified but effective representation of neuronal excitability and firing patterns that is valuable for studying the dynamics of neural networks and related phenomena.

\section{Taylor Approximation to the FHN Model}
\label{Sec:Taylor}
In this section, we describe an algorithm for the solution of the problem defined in (\ref{eq-model1}) with the ICs: $v(0)=0$ and $\omega (0)=-0.2$. This approach gives us an alternative numerical understanding for the problem. Thus, we describe a suitable numerical algorithm for the solution of the problem. With the help of this numerical approach and its convergence, the accurate results are obtained and investigated for the infinity norm ($L_{\infty}$-Norm).

\subsection{Method Description}
We approximate to the system solutions in (\ref{eq-model1}) by truncated Taylor series.
\begin{equation}
\label{Eq:t1}
v(t)\cong v_{N}(t)=\sum_{n=0}^{N}{a_{1,n}(t-c)^n} \quad \mbox{and} \quad \omega(t)\cong \omega_{N}(t)=\sum_{n=0}^{N}{a_{2,n}(t-c)^n},
\end{equation}
where $N$ is the Taylor polynomial degree at the point $t=c$ for $d\leq t\leq e$, for all $d,\, e \in \mathbb{N}$ \cite{wang}. Here we have
\begin{equation}
\label{Eq:t1a}
a_{1,n}=\frac{1}{n!}v^{(n)}(c) \quad \mbox{and} \quad a_{2,n}=\frac{1}{n!}\omega^{(n)}(c), \, n=1,2,\ldots,N,
\end{equation}
the unknown coefficients $a_{1,n}$ and $a_{2,n}$ and we approximate to the numerical solution of $v(t)$ and $\omega(t)$ in the system (\ref{eq-model1}) by using Taylor series approach \cite{bg1,bg2}. Moreover we consider the collocation points defined by\\
\begin{equation}
\label{Eq:t2}
t_i=d+\frac{e-d}{N}i, \quad i=0,1,\ldots, N,
\end{equation}
where we have the pointwise approximation over the interval $[d,e]$, i.e. $d=t_0<t_1<\ldots<t_N=e$. The numerical approach is described by using the matrix relations around $c=0$. Now we consider the matrix forms of the unknown functions, $v(t)$ and $\omega(t)$, and their derivatives \cite{wang,bg3}. Thus the equations in (\ref{eq-model1}) are introduced in the matrix form:\\
\begin{eqnarray}
[v'(t)]&=&\mathbf{T}(t)\mathbf{B}\mathbf{A}_1=\frac{1}{\mu}[f(v,a)]-\frac{1}{\mu}\mathbf{T}(t)\mathbf{A}_2+\mathbf{R},\notag\\
{[\omega'(t)]}&=&\mathbf{T}(t)\mathbf{B}\mathbf{A}_2=\mathbf{T}(t)\mathbf{A}_1-\gamma \mathbf{T}(t)\mathbf{A}_2,  \label{Eq:t3}
\end{eqnarray}
where $\mathbf{R}$ is a constant matrix with respect to $I$. On the other hand, matrix representations of the nonlinear arguments of the system (\ref{eq-model1}) are denoted by $[f(v,a)]$ which will be defined in details. Briefly,
\begin{equation}
\label{Eq:t3a}
[v(t)]=\mathbf{T}(t)\mathbf{A}_{1} \quad \mbox{and} \quad [\omega(t)]=\mathbf{T}(t)\mathbf{A}_{2}.
\end{equation}
Besides, we show the matrix forms of the derivatives of $v(t)$ and $\omega(t)$ in the model (\ref{eq-model1}) \cite{probg1,probg2}.
\begin{equation*}
[v'(t)]=\mathbf{T}(t)\mathbf{B}\mathbf{A}_{1} \quad \mbox{and} \quad [\omega'(t)]=\mathbf{T}(t)\mathbf{B}\mathbf{A}_{2},
\end{equation*}
where
{\small\begin{eqnarray*}
\mathbf{T}(t)&=&\left[\begin{array}{ccccc}
    1 & t& t^2 & \cdots & t^N \\
  \end{array}
  \right], \quad
\footnotesize{ \mathbf{B}=\left[
    \begin{array}{ccccc}
      0 & 1 & 0 & \cdots & 0 \\
      0 & 0 & 2 & \cdots & 0 \\
      \vdots & \vdots & \vdots & \ddots & \vdots \\
      0 & 0 & 0 & \cdots & N \\
      0 & 0 & 0 & \cdots & 0 \\
    \end{array}
  \right]}, \\
\mathbf{A}_{1}&=&\left[
   \begin{array}{ccccc}
     a_{10} & a_{11}& a_{12} & \cdots &a_{1N} \\
   \end{array}
 \right]^{T}, \quad
 \mathbf{A}_{2}=\left[
   \begin{array}{ccccc}
     a_{20} & a_{21}& a_{22} & \cdots &a_{2N} \\
   \end{array}
 \right]^{T}.
\end{eqnarray*}}

On the other hand, we have the function $f(v,a)$ which is defined as
\begin{equation*}
f(v,a)=v(a-v)(v-1)=-v^3+v^2+v^2a-va.
\end{equation*}
Here the nonlinear arguments of the function $f(v,a)$ are also introduced in the matrix forms as follows \cite{cet}:
\begin{equation*}
[f(v,a)]=-\mathbf{T}(t)\overline{\mathbf{T}}(t)\overline{\overline{\mathbf{T}}}(t)\overline{\overline{\mathbf{A}}}_{1}
+a\mathbf{T}(t)\overline{\mathbf{T}}(t)\overline{\mathbf{A}}_{1}-a\mathbf{T}(t)\mathbf{A}_{1}=\mathbf{f}(v,a),
\end{equation*}
where
\begin{equation*}
{[v^{3}]} = \mathbf{T}(t)\overline{\mathbf{T}}(t)\overline{\overline{\mathbf{T}}}(t)\overline{\overline{\mathbf{A}}}_{1}, \quad {[v^{2}]} = \mathbf{T}(t)\overline{\mathbf{T}}(t)\overline{\mathbf{A}}_{1},
\end{equation*}
and
{\small\begin{eqnarray*}
\overline{\mathbf{T}}(t)=
diag\left(\begin{array}{ccccc}
    \mathbf{T}(t), & \mathbf{T}(t),& \cdots, & \mathbf{T}(t) \\
  \end{array}
  \right), \quad
\overline{\overline{\mathbf{T}}}(t)=
diag\left(\begin{array}{ccccc}
   \overline{\mathbf{T}}(t), & \overline{\mathbf{T}}(t),& \cdots, & \overline{\mathbf{T}}(t) \\
  \end{array}
  \right),\\
\overline{\mathbf{A}}_{1}=\left[\begin{array}{ccccc}
    \mathbf{A}_{1} & \mathbf{A}_{1}  & \cdots & \mathbf{A}_{1} \\
  \end{array}
  \right],  \quad
\overline{\overline{\mathbf{A}}}_{1}=\left[\begin{array}{ccccc}
    \mathbf{A}_{1}a_{10} & \mathbf{A}_{1}a_{11}  & \cdots & \mathbf{A}_{1}a_{1N} \\
  \end{array}
  \right].
\end{eqnarray*}}
Hence we have the fundamental matrix relations of the related system in (\ref{eq-model1}). We use the matrix form of the system in (\ref{Eq:t3}) and the collocation points defined in (\ref{Eq:t2}) and we have
\begin{eqnarray}
\mathbf{T}(t_i)\mathbf{B}\mathbf{A}_1&=&\frac{1}{\mu}\mathbf{f}(v,a)-\frac{1}{\mu}\mathbf{T}(t_i)\mathbf{A}_2+\mathbf{R}, \notag\\
\mathbf{T}(t_i)\mathbf{B}\mathbf{A}_2&=&\mathbf{T}(t_i)\mathbf{A}_1-\gamma \mathbf{T}(t_i)\mathbf{A}_2. \label{Eq:t4}
\end{eqnarray}
Alternatively,
\begin{eqnarray}
\left(\mathbf{T}\mathbf{B}+\frac{a}{\mu}\mathbf{T}\right)\mathbf{A}_1
-\frac{a}{\mu}\mathbf{T}\overline{\mathbf{T}}\overline{\mathbf{A}}_{1}
+\frac{a}{\mu}\mathbf{T}\overline{\mathbf{T}}\overline{\overline{\mathbf{T}}}\overline{\overline{\mathbf{A}}}_{1}
-\frac{1}{\mu}\mathbf{T}\mathbf{A}_2&=&\mathbf{R}, \notag\\
-\mathbf{T}\mathbf{A}_1+(\mathbf{T}\mathbf{B}+\gamma \mathbf{T})\mathbf{A}_2&=&\mathbf{0}. \label{Eq:t5}
\end{eqnarray}
Now we write the matrix form of the ICs with the help of the relations in (\ref{Eq:t3a}) and we have
\begin{equation}
\label{Eq:t5a}
[v(0)]=\mathbf{T}(0)\mathbf{A}_{1}=[0], \quad \mbox{and} \quad [\omega(0)]=\mathbf{T}(0)\mathbf{A}_{2}=[-0.2],
\end{equation}
and it is including $[f(v(0),a)]=-\mathbf{T}(0)\overline{\mathbf{T}}(0)\overline{\overline{\mathbf{T}}}(0)\overline{\overline{\mathbf{A}}}_{1}
+a\mathbf{T}(0)\overline{\mathbf{T}}(0)\overline{\mathbf{A}}_{1}-a\mathbf{T}(0)\mathbf{A}_{1}$ for any $a$.
Then we replace the low matrices (\ref{Eq:t5a}) into the alternatively written new matrix system (\ref{Eq:t5}) and we obtain the solutions of the fundamental matrix system in the case the determinant of the system is different than zero. Thus we get the approximate solutions in the form
\begin{equation}
\label{Eq:t6}
v_{N}(t)=\sum_{n=0}^{N}{\frac{1}{n!}v^{(n)}(0)t^n} \quad \mbox{and} \quad \omega_{N}(t)=\sum_{n=0}^{N}{\frac{1}{n!}\omega^{(n)}(0)t^n}.
\end{equation}

\subsection{Convergence Analysis}
Now we perform the error estimation for the numerical solutions of the model defined in (\ref{eq-model1}). Besides, the convergence of the numerical technique based on Taylor truncated series is proved \cite{wang}.

\begin{theorem}
\label{Theo:c1}
  Assume that the function $g(t)$ is defined on the interval $[d,e]$ for the solutions of the systems of the problem defined in (\ref{eq-model1}) together with the ICs, and $u_N(t)$ is the $N$th order Taylor polynomial solution of $u(t)$ which is obtained from (\ref{Eq:t6}) for $n=0,1,...,N$. Hence
  \begin{equation}
  \label{Eq:c1}
    \left\Vert u(t)-u_N(t)\right\Vert_{\infty}\leq\frac{M}{(N+1)!}\left\lvert u^{(N+1)}(\xi)\right\rvert+L\,\max_{0\leq n\leq N} \left\lvert e_n(c) \right\rvert,
  \end{equation}
   where $M=\underset{d\leq t\leq e}{\max}\left \lvert (t-c)^{N+1}\right \rvert$, $L=\norm{\ell}_{\infty}$ and
  $e_n(c)=u^{(n)}(c)-u^{(n)}_{N}(c)$ for $d \leq t\leq e$. In particular, $L_{\infty}$-Norm is defined as $\norm{\ell}_{\infty}=\underset{d\leq t\leq e}{\max}\{\left\lvert l_0(c)\right\rvert, \left\lvert l_1(c)\right\rvert,\ldots,\left\lvert l_n(c)\right\rvert \}$.
\end{theorem}
\begin{proof}
  Let us consider
  \begin{equation*}
    \left\Vert u(t)-u_{N}(t)\right\Vert_{\infty}\leq \left\Vert u(t)-T_n(t) \right\Vert_{\infty} +\left\Vert T_n(t)-u_{N}(t)\right\Vert_{\infty},
  \end{equation*}
  where
  \begin{equation*}
  T_n(t)=\sum_{n=0}^{N}\frac{u^{(n)}(c)(t-c)^n}{n!}, \quad u_N(t)=\sum_{n=0}^{N}\frac{u_N(c)(t-c)^n}{n!}.
  \end{equation*}
Here we consider the residual correction together with the remainder of Taylor polynomial \cite{wang,oliveira}.
\begin{equation*}
  R_n(t)=u(t)-T_n(t)=\frac{u^{(N+1)}(\xi)}{(N+1)!}(t-c)^{N+1}.
\end{equation*}
Thus we have
\begin{equation}
\label{Eq:c2}
  \left\lvert R_n(t)\right\rvert\leq\frac{u^{(N+1)}(\xi)}{(N+1)!}\cdot \underset{d\leq t \leq e}{\max}\left\lvert(t-c)^{N+1}\right\rvert=\frac{M}{(N+1)!}u^{(N+1)}(\xi).
\end{equation}
On the other hand, we have the error function $e_n(c)$ around $t=c$ and the Taylor functions $l_n(c)$ for $n=0, 1,\ldots, N$, respectively:
\begin{equation*}
  e_n(c)=u^{(n)}(c)-u^{(n)}_N(c) \quad \mbox{and}\quad l_n(c)=\frac{(t-c)^n}{n!}.
\end{equation*}
Then we set
\begin{equation*}
  \delta_n=(e_0(c), e_1(c),\ldots, e_n(c), \ldots, e_N(c)), \quad \ell=(l_0(c), l_1(c),\ldots, l_n(c), \ldots, e_N(c))^T.
\end{equation*}
Hence we obtain
\begin{equation}
\label{Eq:c3}
  \left\lvert T_n(t)-u_N(t)\right\rvert=\left\lvert\sum_{n=o}^{N}\left(u^{(n)}(c)-u^{(n)}_N(c)\right)\frac{(t-c)^n}{n!}\right\rvert=\left\lvert \delta_n \cdot \ell\right\rvert=\left\Vert\delta_n\right\Vert_{\infty}\cdot\left\Vert \ell\right\Vert_{\infty}\leq L \left\Vert\delta_n\right\Vert_{\infty}.
\end{equation}
From (\ref{Eq:c2}) and (\ref{Eq:c3}), we have
\begin{eqnarray*}
\left\Vert u(t)-u_N(t)\right\Vert_{\infty}&\leq& \frac{u^{(N+1)}(\xi)}{(N+1)!}\cdot \underset{d \leq t \leq e}{\max}\left\lvert(t-c)^{N+1}\right\rvert+\left\Vert\ell\right\Vert_{\infty}\cdot \left\Vert\delta_n\right\Vert_{\infty}\\
 &=& \frac{M}{(N+1)!}\left\lvert u^{(N+1)}(\xi)\right\rvert+L\, \underset{0 \leq n \leq N}{\max}\left\lvert e_n(c)\right\rvert.
\end{eqnarray*}
This completes the proof.
\end{proof}
Now we have the following Corollary results from (\ref{Theo:c1}):
\begin{corollary}
If $g(t)$ and $u(t)$ run with Taylor polynomial functions, we obtain the Taylor polynomials solutions of (\ref{eq-model1}) together with ICs. Hence we call $u_N(t)$ as the $N$th order Taylor polynomial solution of $u(t)$ for any $N>0$ \cite{wang}.
\end{corollary}
\begin{corollary}
If $g(t)$ and $u(t)$ are $n$th order differentiable functions, then $u_N(t)=T_n(t)$. Namely, we call $u_N(t)$ as the $N$th order Taylor polynomial solution of $u(t)$. Besides, $u_N(t)$ is described as the $N$th order Taylor interpolating polynomial \cite{wang}.
\end{corollary}
Hence, we complete the convergence proof of the numerical method. It is designed based on the approximation to the problem. The analysis of the method introduces therefore a suitable approximation to the problem which is defined for the FHN model. This alternative approach gives a better understanding for the appropriate results which is a novel application for such dynamic model.

\section{Constructed Difference Scheme and its Stability for the FHN Model}
\label{Sec:fdifference}
In this section, we describe an alternative for a numerical scheme and its stability for the solution of the FHN model with the ICs given above. The numerical scheme is based on the finite difference approach. Thus, in this numerical approach we have beneficial outcome for the investigation of the model. The numerical scheme has several steps and we now explain it with its stability. We aim to find approximation to the model and understand it comprehensively. Thus, we introduce an alternative numerical scheme with the stability analysis.
The discretisation of problem (\ref{eq-model1}) is carried out in one step, we
define the grid space \cite{modd}. Now, introduce grids with uniform steps are given as
\begin{equation*}
\overline{W^{\tau }}=\{t_{k}:t_{k}=k\tau ,k=0,1,\ldots ,N,\left. N\tau
=T\right\} ,W^{\tau }=\overline{W^{\tau }}\cap W^{\tau }
\end{equation*}

and introduce the Hilbert space $H=L_{2\tau }=L_{2}(\overline{W^{\tau }})$ of the grid functions $\varphi ^{\tau }(t)=\left\{ \varphi (\tau_{1}i_{1},\tau _{2}i_{2},...,\tau _{k}i_{k}\right\} $ defined on $\overline{%
W^{\tau }}$equipped with the norm $\left\Vert W\right\Vert$, \cite{ashyralyev2015operator}%
\begin{equation}\label{4.1}
\left\Vert \varphi ^{\tau }\right\Vert _{2\tau }=\left( \sum\limits_{t\in
W^{\tau }}\left\vert \varphi ^{\tau }\right\vert \tau _{1}\tau _{2}...\tau
_{k}\right) ^{\frac{1}{2}}.
\end{equation}

Using Taylor expansion for the formula as to $t$ of $v(t)$
and $\omega(t)$, the forward difference formula is obtained as%
\begin{equation}\label{4.2}
\left\{
\begin{array}{l}
v_{t}(t_{k})\approx \frac{%
v_{k+1}-v_{k}}{\tau},\\
\omega_{t}(t_{k})\approx \frac{%
\omega_{k+1}-\omega_{k}}{\tau}.
\end{array}%
\right.
\end{equation}
The ICs can be written as $v_{0}=0,\omega_{0}=-0.2$ for the system (\ref{eq-model1}). The formula (\ref{4.2}) is written into the formula (\ref{eq-model1}), we have%
\begin{equation}\label{4.3}
\left\{
\begin{array}{l}
\omega_{k+1}=(1-\gamma \tau )\omega_{k}+\tau v_{k}, \\
v_{k+1}=v_{k}+\frac{\tau }{\mu }\left[v_{k}(a-v_{k})(v_{k}-1)-\omega_{k}+I\right] \\
v_{0}=0,\omega_{0}=-0.2.%
\end{array}%
\right.
\end{equation}

\begin{figure}[H]
		\centering
		\includegraphics[width=1\textwidth]{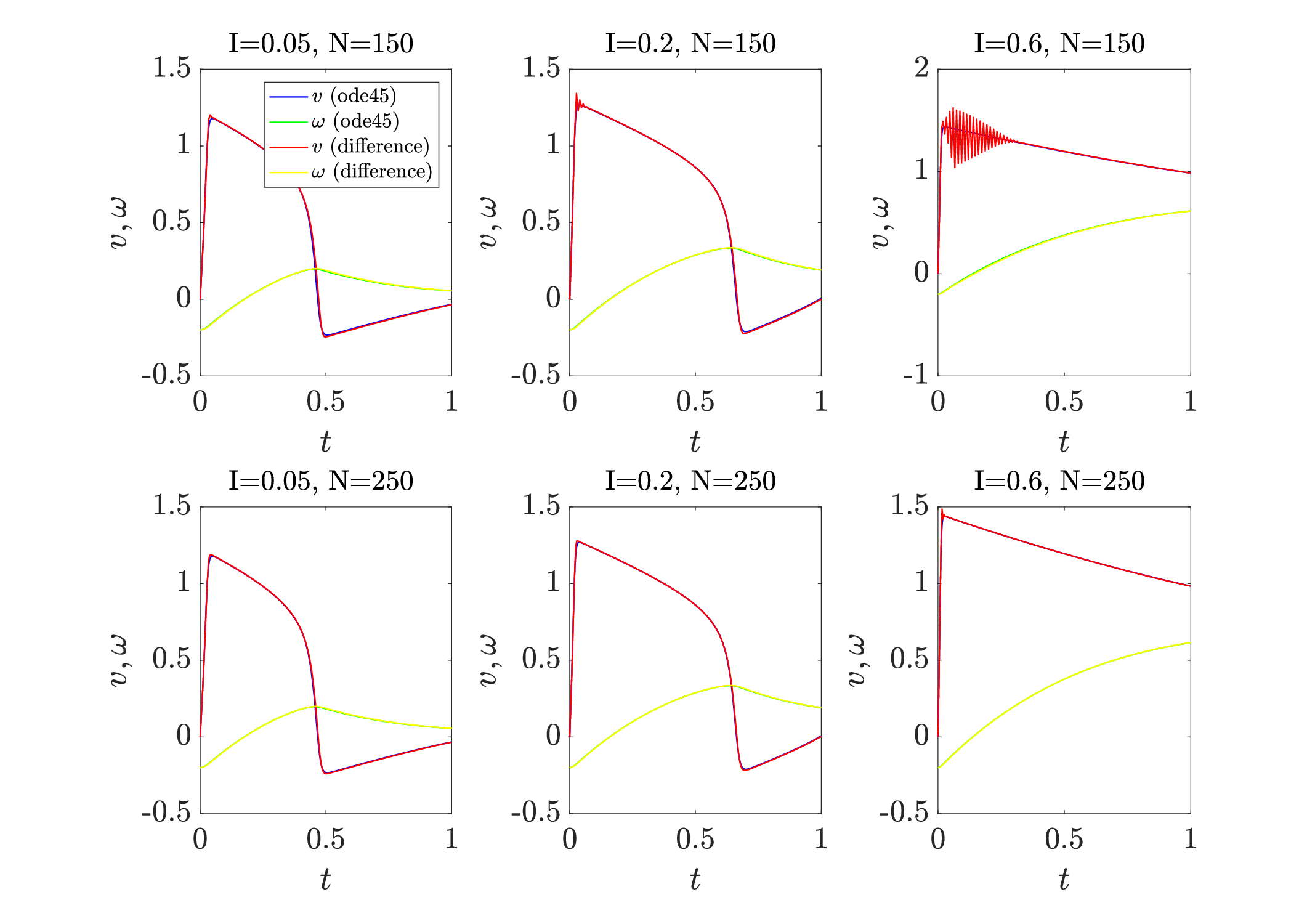}
		\caption{Gives the approximation solutions of $v(t)$ and $\omega(t)$ for $\mu=0.008$, $\gamma=1.18$ , $a=0.22$, and  $0\le t \le1$ .}
		\label{fig22}
\end{figure}

\begin{figure}[H]
		\centering
		\includegraphics[width=1\textwidth]{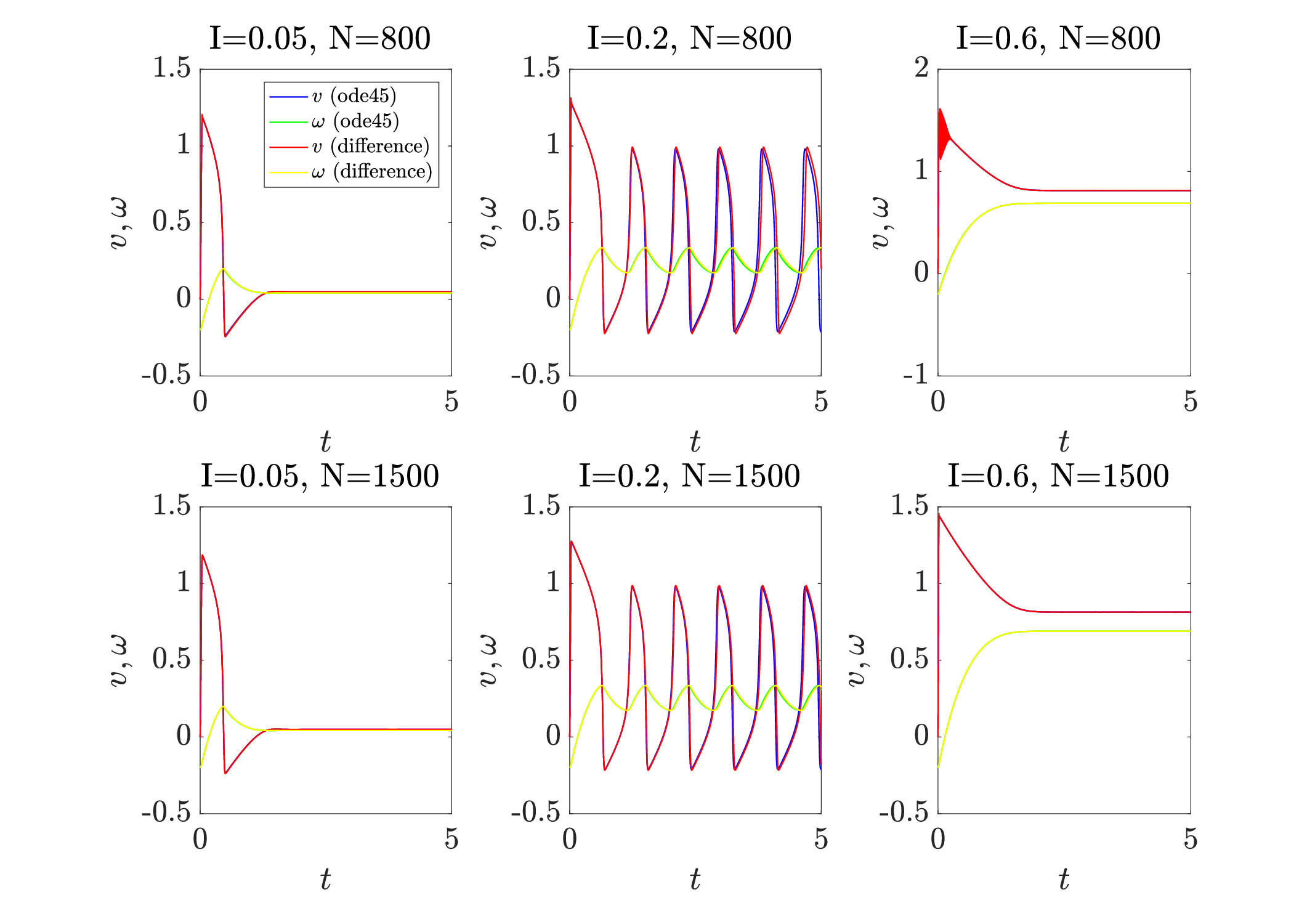}
		\caption{Gives the approximation solutions of $v(t)$ and $\omega(t)$ for $\mu=0.008$, $\gamma=1.18$ , $a=0.22$, and  $0\le t \le 5$ .}
		\label{fig33}
\end{figure}

Now, we will prove the theorem of stability estimates for the difference scheme formula (\ref{4.3}).\\
\begin{theorem}
The following stability estimates is satisfies for the formula (\ref{4.3})
with depend on ICs as:
\begin{enumerate}[label=(\roman*)]
\item $\left\Vert v_{k}\right\Vert _{H}\leq \frac{\tau }{\mu }%
\left[\sum\limits_{j=1}^{k}(\left\Vert f_{j-1}-\omega_{j-1}\right\Vert _{H})+k\left\Vert I\right\Vert _{H}\right]$
\item $\left\Vert \omega_{k}\right\Vert _{H}\leq \tau
\sum\limits_{j=1}^{k-1}\left\Vert v_{j}\right\Vert _{H},$
\end{enumerate}
\end{theorem}
\begin{proof}
\begin{enumerate}[label=(\roman*)]
\item From the formula (\ref{4.3}), we can write%
\begin{equation*}
v_{k+1}=v_{k}+\frac{\tau }{\mu }[f_{k}-\omega_{k}+I],
\end{equation*}%
and form that%
\begin{equation}\label{4.4}
v_{k}=v_{k-1}+\frac{\tau }{\mu }[f_{k-1}-\omega_{k-1}+I],
\end{equation}%
here $f_{k}=v_{k}(v_{k}-a)(v_{k}-1).$ The formula (\ref{4.4}) is clearly written as%
\begin{eqnarray}\label{4.5}
v_{1} &=&v_{0}+\frac{\tau }{\mu }[f_{0}-\omega_{0}+I],\text{ }v_{0}=0;  \notag \\
v_{2} &=&\frac{\tau }{\mu }[f_{0}-\omega_{0}+I]+\frac{\tau }{\mu }[f_{1}-\omega_{1}+I],
\notag \\
&&...  \notag \\
v_{k} &=&\frac{\tau }{\mu }[f_{0}+f_{1}+...+f_{k-1}-(\omega_{0}+\omega_{1}+...+\omega_{k-1})+I+I+...+I],  \notag \\
v_{k} &=&\frac{\tau }{\mu }[\sum\limits_{j=1}^{k}(f_{j-1}-\omega_{j-1})+kI].
\end{eqnarray}%
Using triangle inequality and applying the norm (\ref{4.1}) to the formula (%
\ref{4.5}), we can obtain%
\begin{equation}\label{4.6}
\left\Vert v_{k}\right\Vert _{H}\leq \frac{\tau }{\mu }
[\sum\limits_{j=1}^{k}(\left\Vert f_{j-1}-\omega_{j-1}\right\Vert _{H})+k\left\Vert I\right\Vert _{H}].
\end{equation}%
From the extreme points of a function, it can be easily seen that the
maximum point of $f_{k}$ is $\frac{2}{3}$and the minimum point of $\omega_{k}$ is
$-\frac{1}{5}.$ We know from the first and second values that $f_{k}$ and $\omega_{k}$ functions are increasing functions. Thus, we have%
\begin{equation}\label{4.7}
\max_{1\leq j\leq k}\left\Vert f_{j-1}-\omega_{j-1}\right\Vert _{H}=\frac{2}{3}+%
\frac{1}{5}=\frac{13}{15}.
\end{equation}%
Considering that $1<k<N,$ $\tau =\frac{1}{N}$ and $\mu $ is finite value, we obtain%
\begin{equation}\label{4.8}
\left\Vert v_{k}\right\Vert _{H}\leq (\frac{13\tau }{15\mu }+k\left\Vert I\right\Vert _{H}).
\end{equation}%
From (\ref{4.8}), the stability estimates are satisfied for (\ref{4.8}).

\item For the proof of the theorem, the formula (\ref{4.2}) can be rewritten
as%
\begin{eqnarray}
\label{4.9}
\omega_{0}&=&-\frac{1}{5}, \\
\omega_{1}&=&-\frac{1}{5}(1-\gamma \tau ), \\
\omega_{2}&=&-\frac{1}{5}(1-\gamma \tau )^{2}+v_{1}, \\
\omega_{3}&=&-\frac{1}{5}(1-\gamma \tau )^{3}+\tau (1-\gamma \tau )v_{1}+\tau v_{2}, \\
\omega_{k}&=&-\frac{1}{5}(1-\gamma \tau )^{k}+\tau (1-\gamma \tau )^{k-2}v_{1}+\tau (1-\gamma \tau )^{k-3}v_{2}+...+v_{k-1}\tau .%
\end{eqnarray}%
Using triangle inequality and applying the norm (\ref{4.1}) to the formula (\ref{4.9}), we can obtain%
\begin{eqnarray}\label{4.10}
\left\Vert \omega_{k}\right\Vert _{H} &=&\left\vert -\frac{1}{5}(1-\gamma \tau
)^{k}\right\vert +\left\vert \tau (1-\gamma \tau )^{k-2}\right\vert
\left\Vert v_{1}\right\Vert _{H}+...  \notag \\
&&+\left\vert \tau (1-\gamma \tau )^{k-3}\right\vert \left\Vert
v_{2}\right\Vert _{H}+...+\tau \left\Vert v_{k-1}\right\Vert _{H}  \notag \\
&\leq &1+\tau \lbrack \left\Vert v_{1}\right\Vert _{H}+\left\Vert
v_{2}\right\Vert _{H}+...+\left\Vert v_{k-1}\right\Vert _{H}]  \notag \\
&\leq &\tau \sum\limits_{j=1}^{k-1}\left\Vert v_{j}\right\Vert _{H}.
\end{eqnarray}%
Using the formula (\ref{4.8}), the stability estimates are satisfied for the
formula (\ref{4.10}). Thus, we complete the proof of the theorem.
\end{enumerate}
\end{proof}
Let us now present the numerical results of the problem defined in (\ref{eq-model1}) together with the ICs.
By using the formula (\ref{4.3}) and MATLAB programming, we have the simulations which are shown in Figures (\ref{fig22}) and (\ref{fig33}), respectively.

%
%
%
%
%
%
From Figures \ref{fig22} and \ref{fig33}, $\omega(t)$ increases up to given value in the formula (\ref{4.8}), while $v(t)$ decreases up to that value. After this value, it continues in parallel as depending on the formula (\ref{4.10}). However, as the number of $N$ intervals increases, the physical appearance of the figures becomes clearer as the $\omega$ and $v$ functions are closer to their real values. At the same time, as the $t$ time interval widens, the behavior of the figures becomes more pronounced, for example from \ref{fig22} to \ref{fig33}.

We now investigate the convergency rates of the numerical approach in (\ref{Sec:Taylor}) for the $L_{\infty}$-Norm \cite{feng,avaz}. This gives us a better understanding for the Taylor approach together with the Finite Difference approximation results. Here we also consider the parameters as follows: $I=0.6$, $\mu=0.008$, $\gamma=1.18$, $a=0.22$, and $\tau=0.00025$, for  $0\le t \le1$. From Table (\ref{tab1}) and Table (\ref{tab2}), we can easily see that the numerical approach ensures efficient results for higher order $N$ values. Thus the numerical approach has beneficial findings whenever we reach for larger iterations. We can also understand that the approximations for the different truncation limits have analogue results at the $L_{\infty}$-Norm.

\begin{table}[H]
\centering
\caption{Error convergence of approximated $v(t)$ for the $L_{\infty}$-Norm and $N=4$, $5$ and $6$, respectively.}
\label{tab1}
\begin{tabular}{cccc} \hline
\textbf{$t$} & \textbf{$N=4$} & \textbf{$N=5$} & \textbf{$N=6$}\\ \hline
0.0 & 0.10620E-3 & 0.11630E-4 & 0.13978E-5 \\
0.1 & 0.11560E-3 & 0.15122E-4 & 0.13596E-5 \\
0.2 & 0.13323E-3 & 0.18851E-4 & 0.13671E-5 \\
0.3 & 0.10003E-3 & 0.14520E-4 & 0.14582E-5 \\
0.4 & 0.11588E-3 & 0.15505E-4 & 0.14896E-5 \\
0.5 & 0.13994E-3 & 0.11091E-4 & 0.19071E-5 \\
0.6 & 0.14278E-3 & 0.19985E-4 & 0.13601E-5 \\
0.7 & 0.13775E-3 & 0.20120E-4 & 0.14772E-5 \\
0.8 & 0.15953E-3 & 0.20512E-4 & 0.14850E-5 \\
0.9 & 0.15029E-3 & 0.20001E-4 & 0.12417E-5 \\
1.0 & 0.19003E-3 & 0.35277E-4 & 0.11230E-4 \\ \hline
\end{tabular}
\end{table}

\begin{figure}[ht!]
\centering
\begin{tabular}{cc}
\includegraphics[width=0.45\textwidth]{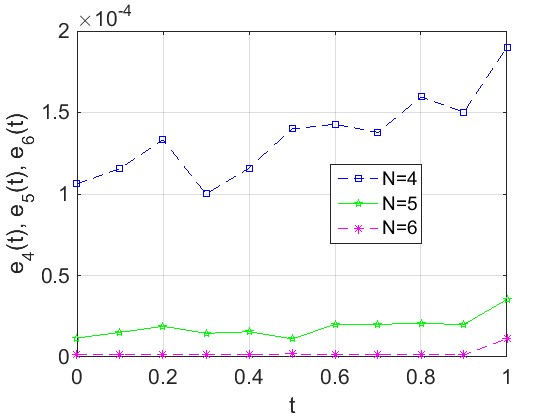} \\
\textbf{(a)}  \\[6pt]
\end{tabular}
\begin{tabular}{cc}
\includegraphics[width=0.45\textwidth]{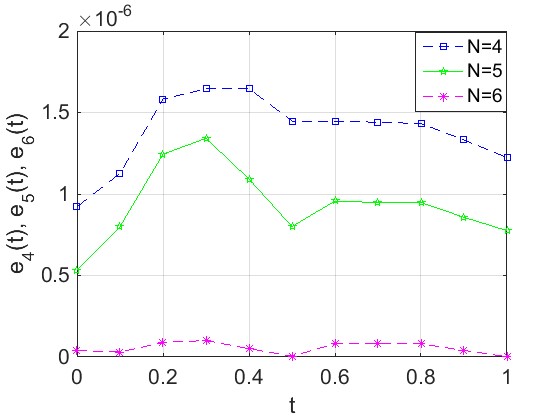} \\
\textbf{(b)}\\[6pt]
\end{tabular}
\caption{The $L_{\infty}$-Norm of error functions for $N=4$, $5$ and $6$, respectively, where(a) shows the results for $v(t)$ and (b) $\omega(t)$.}
\label{Fig:Errors}
\end{figure}

Besides Figure \ref{Fig:Errors} shows us that the behaviour of the convergency results shows similarities even for the different unknowns in the system, $v(t)$ and $\omega(t)$. On the other hand, these $L_{\infty}$-Norm of error functions show that the behaviour of these functions vary for the same interval where we have $[0,1]$.

In particular, the $L_{\infty}$-Norm of error function of $v(t)$ has complication which has been also basically seen in Figure \ref{Fig:Errors}. Namely, the unknown function has the $L_{\infty}$-Norm of its error function has increasing outcome which is visible in Figure \ref{Fig:Errors} and at Table \ref{tab1}. There is also seen that $\omega(t)$ has been simply exponentially growing function. Thus the $L_{\infty}$-Norm of error function of $\omega(t)$ has decreasing error results which is seen in Figure \ref{Fig:Errors} and at Table \ref{tab2}. The computation time for various numerical schemes has been analysed for \(N\) and simulation time \(t\) \cite{olmos09}. The calculation per unit (CPU) (s) is obtained for the comparison values for the speed of $t$. As an example, for the fixed parameters, ICs and $N=4,5$ and $6$ for $\omega(t)$ at $0\leq t \leq 1$, we obtain the CPU (s): $0.111$ s, $0.175$s,  and $0.188$s, respectively.

\begin{table}[H]
\centering
\caption{Error convergence of approximated $\omega(t)$ for the $L_{\infty}$-Norm and $N=4$, $5$ and $6$, respectively.}
\label{tab2}
\begin{tabular}{cccc} \hline
\textbf{$t$} & \textbf{$N=4$} & \textbf{$N=5$} & \textbf{$N=6$}\\ \hline
0.0 & 0.09217E-5 & 0.05309E-5 & 0.04012E-6 \\
0.1 & 0.11259E-5 & 0.08000E-5 & 0.03015E-6 \\
0.2 & 0.15822E-5 & 0.12443E-5 & 0.08928E-6 \\
0.3 & 0.16448E-5 & 0.13418E-5 & 0.10092E-6 \\
0.4 & 0.16449E-5 & 0.10919E-5 & 0.05097E-6 \\
0.5 & 0.14476E-5 & 0.80021E-6 & 0.04902E-7 \\
0.6 & 0.14445E-5 & 0.95708E-6 & 0.81000E-7 \\
0.7 & 0.14409E-5 & 0.94712E-6 & 0.85011E-7 \\
0.8 & 0.14319E-5 & 0.94810E-6 & 0.80031E-7 \\
0.9 & 0.13318E-5 & 0.85701E-6 & 0.40122E-7 \\
1.0 & 0.12231E-5 & 0.77500E-6 & 0.01884E-7 \\ \hline
\end{tabular}
\end{table}

\section{Conclusion}\label{Sec:Conc}
The study purposes a comprehensive analysis of the FHN model problem defined in (\ref{eq-model1}) with ICs. Firstly, the model is revisited with its equilibrium and stability analysis via a characteristic polynomial. In fact, the linearisation of this nonlinear system near the equilibrium points is key to determine the stability with the eigenvalues found by the Jacobian matrix. Depending on the value of the external input current ($I$), stable and unstable dynamics through a limit cycle around the positive equilibrium are shown as an example. 
The model results have been studied and confirmed by the authors, see for example \cite{izhikevich2006fitzhugh,faghih2010fitzhugh,kostova2004fitzhugh}.

The highlight of this manuscript is to introduce a numerical method based on the truncated Taylor series and provide alternative approaches to obtain the solutions and stability of the system. In this regard, Taylor polynomial approach has been applied on the model problem and we have the convergency result of the method. Moreover, finite difference scheme has been performed and numerical simulations of the results have been obtained. Thus, the FHN model is motivated by the numerical approaches which show us appropriate dynamical results of the unknown functions in the model. The straightforwardness of the combined technical application gives a novelty for the scientific computing approach. Namely, we have completed a comprehensive numerical application with a highly motivated FHN model described in details with its stability analysis and the bifurcation approach. This consequence is supported by the figures and the tables where we can easily see the beneficial outcome. The results can be also enriched by the different values of the parameters and adapted for the different systems. The numerical scheme based on Taylor polynomials is of advantageous outcomes regarding the convergency results and the explicit procedure in the concept of an algorithmic approach \cite{gurbuz2021rumour}.
Although this paper presents a thorough investigation of the FHN model using analytical and numerical techniques, certain limitations should be noted. To account for complicated dynamics, approximations are used in the linear stability analysis and local stability conditions. Precision and computational resource limitations affect the numerical simulations. Although developed for this model, the innovative Taylor polynomial and built difference scheme approaches have not been as effective or scalable to larger frameworks. The FHN model is also a simplistic representation that may not fully capture all biological nuances. Despite these drawbacks, the work provides a useful framework for studying relaxation oscillators and opens the door for further studies that address these limitations, consider different models and improve numerical techniques. On the other hand, finite difference scheme for the problem concluded the approximation with highly reasonable results which is uniquely designed as an alternative approach for the problem. The combination of the numerical schemes for the FHN model open new tasks for the future investigation and introduce a novel approximation for the which will be essential for the further research on approximation for such dynamical system problems.

Furthermore, the model can be dealt with its different aspects. In particular, the numerical scheme applied here can be considered for the stochastic analysis of the FHN model as a future Outlook. Besides, the simplification of the FHN model can be obtained by using a piecewise-linear transformation. Thus we can get mathematically more attainable McKean model and we compare the results from the theoretical and numerical aspects \cite{mckean1970nagumo}. Following the ideas presented in \cite{gokcce2021mathematical,gurbuz2021analysis}, the extension of the study is also possible with a spatial diffusion term which is to describe the voltage variable and we compare the dynamics in terms of various theoretical and numerical schemes of partial differential equations and travelling wave analysis by including the time delay.
In addition, by adapting the model and using data-driven methods, such as any deep learning approach, it can also be applied to in-depth investigations that are in high demand \cite{royal22,kossa23}.

\bibliographystyle{unsrtnat}
\bibliography{fhnbg.bst}  






\end{document}